\documentclass[12pt]{amsart}
\usepackage[colorlinks=true, pdfstartview=FitV, linkcolor=blue, citecolor=blue, urlcolor=blue, breaklinks=true]{hyperref}
\usepackage{amsmath,amsfonts,amssymb,amsthm,amscd,comment}
\usepackage{latexsym}
\usepackage[usenames]{color}
\usepackage[all]{xy}
\usepackage{tikz}
\newcommand\C{\mathbb{C}}
\newcommand\Z{\mathbb{Z}}

\newcommand\R{\mathbb{R}}
\newcommand\F{\mathbb{F}}
\newcommand\N{\mathbb{N}}

\newcommand\id{\mathrm{id}}

\newcommand\Sym{\mathrm{Sym}}

\newcommand\rmH{\mathrm{H}}
\newcommand\bimod{\mathrm{bimod}}
\newcommand\Fun{\mathrm{Fun}}
\newcommand\GL{\mathrm{GL}}
\newcommand\SL{\mathrm{SL}}

\newcommand\fh{\mathfrak{h}}

\newcommand\bF{\mathbf{F}}

\newcommand\bA{\mathbf{A}}

\newcommand\cH{\mathcal{H}}
\newcommand\cF{\mathcal{F}}

\newcommand{\G}{\Gamma}

\DeclareMathOperator{\Hom}{Hom}

\DeclareMathOperator{\End}{End}

\DeclareMathOperator{\Ind}{Ind}
\DeclareMathOperator{\Res}{Res}
\DeclareMathOperator{\tr}{tr}

\theoremstyle{plain}
\newtheorem{theo}{Theorem}[section]
\newtheorem{prop}[theo]{Proposition}

\theoremstyle{definition}
\newtheorem{defin}[theo]{Definition}
\newtheorem*{rem*}{Remark}

\numberwithin{equation}{section}

\allowdisplaybreaks

\begin{document}
%
%%%%%%%%%%%%%%%%%%%%%%%%%%%%%%%%%%%

\title[A survey of Heisenberg categorification via graphical calculus]{A survey of Heisenberg categorification \\ via graphical calculus}

\author{Anthony Licata}
\address{A.~Licata: Department of Mathematics, Stanford University}

\author{Alistair Savage}
\address{A.~Savage: Department of Mathematics and Statistics, University of Ottawa}
\thanks{The research of the second author was supported by a Discovery Grant from the Natural Sciences and Engineering Research Council of Canada.}

\subjclass[2010]{20C08, 14F05, 17B65}

%\date{\today}

\begin{abstract}
In this expository paper we present an overview of various graphical categorifications of the Heisenberg algebra and its Fock space representation.  We begin with a discussion of ``weak'' categorifications via modules for Hecke algebras and ``geometrizations'' in terms of the cohomology of the Hilbert scheme of points on the resolution of a simple singularity.  We then turn our attention to more recent ``strong'' categorifications involving planar diagrammatics and derived categories of coherent sheaves on Hilbert schemes.
\end{abstract}

\maketitle \thispagestyle{empty}

\tableofcontents

%%%%%%%%%%%%%%%%%%%%%%%%%%%%%%%%%%%%%%%%%%%%%%%%%%%%%%%%%%%%%%%%%%%
\section*{Introduction}
%%%%%%%%%%%%%%%%%%%%%%%%%%%%%%%%%%%%%%%%%%%%%%%%%%%%%%%%%%%%%%%%%%%

Heisenberg algebras play a fundamental role in quantum physics.  The \emph{Heisenberg algebra of rank $r$} is the unital associative $\C$-algebra with generators $p_n, q_n$, $1 \le n \le r$, and relations
\begin{equation} \label{eq:CCR}
  p_nq_m = q_mp_n + \delta_{n,m}1, \quad p_np_m = p_mp_n, \quad	q_nq_m = q_mq_n,\quad 1 \le n,m \le r.
\end{equation}
In the physics literature, the relations~\eqref{eq:CCR} are sometimes called the \emph{canonical commutation relations} (often with different constants).  Physically, the generators $p_n$ and $q_n$ correspond (up to scalar multiples) to position and momentum operators in a single particle system with $r$ degrees of freedom.  They are also crucial in the study of the quantum harmonic oscillator, the quantum-mechanical analogue of the classical harmonic oscillator.  The physical importance of this system comes from the fact that it can be used to approximate an arbitrary potential in a neighborhood of an equilibrium point, together with the fact that it is one of the few quantum-mechanical systems for which a simple, exact solution is known.

The Heisenberg algebra has an irreducible representation, called \emph{Fock space}, which plays a prominent role in quantum mechanics.  The Stone-von Neumann Theorem asserts that Fock space is the unique irreducible representation of the Heisenberg algebra generated by \emph{vacuum vectors} annihilated by the $p_n$.  This important theorem was a key step in the early understanding of quantum mechanics in that it showed that the Schr\"odinger wave formulation of quantum mechanics and the Heisenberg matrix formulation are physically equivalent.  We refer the reader to the expository article \cite{Ros04} for further details on the history of this theorem and its implications for physics.

In the current paper, we are concerned with an infinite rank version of the Heisenberg algebra with generators $p_n,q_n$, $n \in \N_+ := \Z_{>0}$, and relations
\[
  p_nq_m = q_mp_n + \delta_{n,m}1, \quad p_np_m = p_mp_n, \quad	q_nq_m = q_mq_n, \quad n,m \in \N_+.
\]
This algebra plays a key role in quantum field theory and in the representation theory of infinite-dimensional Lie algebras.  Our aim is to give an overview of some of the categorifications of this important algebra.

In general, categorification is a process in which sets are replaced by categories and equalities by isomorphisms.  In categorifying an algebraic structure, one seeks to find a category that recovers it after passing to the Grothendieck group.  The original algebraic object is then seen to be a shadow of higher categorical structure.  In recent years, there has been considerable interest in the categorification of knot invariants, representations of Lie algebras, and quantum groups.

We begin in Section~\ref{sec:heisenberg-defs} by introducing various incarnations of the Heisenberg algebra and the Fock space, including quantizations and integral versions of both.  Then, in Section~\ref{sec:weak-cat-and-geom}, we discuss some of the first indications of the existence of interesting categorifications of the Heisenberg algebra.  In particular, we describe the work of Geissinger and Zelevnisky, who realized the Fock space as the Grothendieck group of the category of modules for symmetric groups and Hecke algebras.  As we explain, these are examples of ``weak'' categorifications.  We also briefly recall in Section~\ref{sec:weak-cat-and-geom} the realizations, due to Nakajima and Grojnowski, of the Fock space in terms of the cohomology of the Hilbert scheme of points on the resoluation of a simple singularity of type ADE.  This gives a geometric realization, or ``geometrization'', of the Fock space.
As we will see, these weak categorifications and geometrizations can be lifted to ``strong'' categorifications.

In Section~\ref{sec:graphical-cat} we begin our treatment of strong categorifications of the Heisenberg algebra.  We first describe a recent categorification in terms of planar diagrammatics, which $q$-deforms a construction of Khovanov \cite{Kho10}.  This can be seen as the strong analogue of the aforementioned weak categorification via modules for Hecke algebras.  In Section~\ref{sec:quantum-heis-cat} we present a strong categorification, also via planar diagrammatics, which arises from the geometrizations mentioned above.  Both the categorifications from Section~\ref{sec:graphical-cat} and Section~\ref{sec:quantum-heis-cat} were inspired by the work of Khovanov, who initiated the study of graphical categorification of the Heisenberg algebra in \cite{Kho10}.

There is much work on the geometrization and categorification of the Heisenberg algebra and its Fock space that we are not able to cover in detail in the current paper.  We therefore conclude in Section~\ref{sec:further} with a brief overview of some other work appearing in the literature.

\subsection*{Acknowledgements}

The second author would like to thank the Institut de Math\'ematiques de Jussieu and the D\'epartement de Math\'ematiques d'Orsay for their hospitality during his stays there, when the writing of the current paper took place.

%%%%%%%%%%%%%%%%%%%%%%%%%%%%%%%%%%%%%%%%%%%%%%%%%%%%%%%%%%%%%%%%%%%
%
\section{The Heisenberg algebra and Fock Space}
\label{sec:heisenberg-defs}
%
%%%%%%%%%%%%%%%%%%%%%%%%%%%%%%%%%%%%%%%%%%%%%%%%%%%%%%%%%%%%%%%%%%%

%%%%%%%%%%%%%%%%%%%%%%%%%%%%%%%%%%%
\subsection{The Heisenberg algebra}
%%%%%%%%%%%%%%%%%%%%%%%%%%%%%%%%%%%

The \emph{Heisenberg algebra} $\fh$ in infinitely many variables is the unital associative $\C$-algebra with generators $p_n,q_n$, $n \in \N_+$, and relations
\begin{equation} \label{eq:usual-h-relations}
  p_nq_m = q_mp_n + \delta_{n,m}1, \quad p_np_m = p_mp_n, \quad	q_nq_m = q_mq_n, \quad n,m \in \N_+.
\end{equation}
The algebra $\fh$ occurs naturally in mathematics in several different variations and with different presentations.  We recall some of them here.  First, occasionally the generators are rescaled in such a way that the presentation is in terms of generators $p_n$, $n \in \Z \setminus \{0\}$, and relations
\begin{equation} \label{eq:usual-h-relations'}
  p_np_m = p_mp_n + n\delta_{n,-m}1, \quad n,m \in \Z \setminus \{0\}.
\end{equation}

In another alternative presentation of the algebra $\fh$, less obviously equivalent to the presentation  \eqref{eq:usual-h-relations}, the generators are $a_n,b_n$, $n \in \N_+$, and the relations are
\begin{equation} \label{eq:modified-h-relations}
	a_nb_m	=	b_ma_n + b_{m-1}a_{n-1}, \quad a_na_m = a_ma_n,\quad	b_nb_m = b_mb_n, \quad n,m \in \N_+,
\end{equation}
(see \cite[\S1]{Kho10} and \cite[\S1]{LS11}).  In the above, we declare $b_0 = a_0 = 1$.  The unital ring $\fh_\Z \subseteq \fh$ generated by the $a_n$ and $b_n$ is an \emph{integral form} of the Heisenberg algebra, i.e.\ $\C \otimes_\Z \fh_\Z \cong \fh$.

%%%%%%%%%%%%%%%%%%%%%%%%%%%%%%%%%%%%%%%%
\subsection{Lattice Heisenberg algebras} \label{sec:lattice-heisenberg}
%%%%%%%%%%%%%%%%%%%%%%%%%%%%%%%%%%%%%%%%

Let $L$ be a lattice, that is, a finite rank free abelian group equipped with a symmetric bilinear form $\langle \cdot,\cdot \rangle: L\times L \longrightarrow \Z$.  Fix a basis $\alpha_1,\hdots,\alpha_k$ of $L$.  The \emph{lattice Heisenberg algebra} $\fh^L$ associated to $L$ is then defined to be the unital algebra with generators $p_{i,n}$, $i \in \{1,\dots,k\}$, $n\in \Z \setminus \{0\}$, and relations
\begin{equation} \label{eq:usual-h-relations-lattice}
  p_{i,n}p_{j,m} = p_{j,m}p_{i,n} + n\delta_{n,-m}\langle \alpha_i,\alpha_j\rangle 1, \quad i,j \in \{1, \cdots, k\},\ n,m \in \Z \setminus \{0\}.
\end{equation}
Moreover, for any $v=\sum_i m_i\alpha_i\in L$, we may define the element $p_{v,n} \in \fh^L$ by linearity: $p_{v,n} = \sum m_i p_{i,n}$.  In particular, the isomorphism class of the algebra $\fh^L$ does not depend on the choice of basis of $L$.
When $L=\Z$ and the bilinear form is multiplication, this definition agrees with the definition of $\fh$ given in \eqref{eq:usual-h-relations'}.  Heisenberg algebras associated with different lattices show up naturally in a variety of contexts.  An important specific example of a lattice Heisenberg algebra comes from the case when the lattice is associated to a simply-laced finite or affine Dynkin diagram.  More specifically, let $I_\G$ denote the set of nodes of a simply-laced Dynkin diagram $\Gamma$ of finite or affine type (recall that the diagram in affine type has one more node than the diagram of the corresponding finite type).  We let $L_\G = \Z^{I_\G}$ be the free $\Z$-module spanned by $I_\G$. To simplify notation, we denote the basis element corresponding to $i\in I_\G$ by $i$ (as opposed to $\alpha_i$).  We equip $L_\G$ with a symmetric bilinear form $\langle \cdot,\cdot \rangle$ by defining
\[
	\langle i,j \rangle=
	\begin{cases}
       2 & \text{ if } i=j, \\
        -1 & \text{ if } i \ne j \text{ are connected by an edge,} \\
	0 & \text{ if } i \ne j \text{ are not connected by an edge, }
\end{cases}
\]
and extending to all of $L_\G$ bilinearly.  Thus the matrix $(\langle i,j\rangle)_{i,j}$ is the corresponding finite or affine Cartan matrix.  When $\G$ is an affine Dynkin diagram, the associated Heisenberg algebra $\fh^{L_\G}$ is sometimes called a \emph{toroidal Heisenberg algebra}.  These algebras play an important role in the representation theory of infinite-dimensional Lie algebras and mathematical physics.  In particular, one important feature of the Heisenberg algebra $\fh^{L_\G}$ is that it admits a quantization, that is, a deformation over $\C[t,t^{-1}]$.  The quantum Heisenberg algebra $\fh^{L_{\G,t}}$ is defined to be the unital algebra generated by $p_{i,n}$, $i \in \{1,\dots,k\}$, $n\in \Z \setminus \{0\}$, with relations
\[
	p_{i,n}p_{j,m} = p_{j,m}p_{i,n} + \delta_{n,-m}[n\langle i,j \rangle] \frac{[n]}{n} 1.
\]
Here $[k+1] = t^{-k} + t^{-k+2} + \hdots + t^{k-2} + t^k$ denotes the quantum integer.  Note that when $t=1$, the quantum Heisenberg algebra $\fh^{L_\G,t}$ becomes the ordinary Heisenberg algebra
$\fh^{L_\G}$.  Some literature changes the above relations by introducting a minus sign in front of the term $\delta_{n,-m}[n\langle i,j \rangle] \frac{[n]}{n} 1$, though this change does not change the isomorphism class of the resulting algebra.

For $n\geq 0$ and $i\in I_\G$, we define new elements $a_i^{(n)},b_i^{(n)}\in \fh^{L_\G,t}$ using  the generating functions
\[
	\exp\big(\sum_{m\geq 1} \frac{p_{i,-m}}{[m]} z^m) = \sum_{n\geq 0} b_i^{(n)} z^n
	\text{ and }
	\exp\big(\sum_{m\geq 1} \frac{p_{i,m}}{[m]} z^m) = \sum_{n\geq 0} a_i^{(n)} z^n.
\]
The elements $\{a_i^{(n)},b_i^{(n)}\}$ also generate $\fh^{L_\G,t}$.  As shown in \cite{CauLic10}, the defining relations in $\fh^{L_\G,t}$ for these new generators are
\begin{align*}
    a_i^{(n)}   a_j^{(m)}  &= a_j^{(m)}  a_i^{(n)} , \\
    b_i^{(n)}  b_j^{(m)} &= b_j^{(m)} b_i^{(n)}, \\
   b_i^{(n)} a_i^{(m)} &= \sum_{k\geq 0} [k+1] a_i^{(m-k)}b_i^{(n-k)},\\
  b_i^{(n)} a_j^{(m)} &= \sum_{k=0,1} a_j^{(m-k)}b_i^{(n-k)}, \text{ when }
  \langle i,j\rangle = -1,\\
  b_i^{(n)} a_j^{(m)} &= a_j^{(m)}b_i^{(n)}  \text{ when }   \langle i,j\rangle = 0.
\end{align*}
By convention, we set $a_i^{(n)} = b_i^{(n)} = 0$ for $n<0$.  Thus the summations in the above relations are finite.

We note that much of the representation theory literature on infinite-dimensional Lie algebras considers the Heisenberg Lie algebra rather than the Heisenberg algebra $\fh$, which is the enveloping algebra of the Lie algebra.  Moreover, the unit $1$ in the above relations is often replaced by a central generator $c$.

%%%%%%%%%%%%%%%%%%%%%%%%%%%%%%%%%%%%%%%%%%
\subsection{The Fock Space representation}
%%%%%%%%%%%%%%%%%%%%%%%%%%%%%%%%%%%%%%%%%%

Let $\fh^-$ denote the unital subalgebra of $\fh$ generated by the $p_{-i}$, $i \in \N_+$, and let $\fh^+$ denote the subalgebra of $\fh$ generated by the $p_i$, $i \in \N_+$.  We then have an isomorphism of vector spaces $\fh \cong \fh^+\otimes \fh^-$.  Let $W$ denote the trivial representation of $\fh^+$.  The induced representation
\[
	\cF = \fh\otimes_{\fh^+} W
\]
is known as the \emph{Fock space representation} of the Heisenberg algebra.  As a vector space, $\cF$ is isomorphic to $\fh^-$, which is a polynomial algebra in the generators $p_{-i}$, $i \in \N_+$:
\[
	\cF \cong \C[p_{-1},p_{-2},\dots].
\]
Moreover, the representation $\cF$ is faithful and the Heisenberg algebra $\fh$ may be thought of as an algebra of differential operators on the above polynomial algebra.  In this realization, $p_{-i} \in \fh$, $i \in \N_+$, corresponds to multiplication by the variable $p_{-i}$, and the $p_i$ of presentation~\eqref{eq:usual-h-relations'} corresponds to $i \frac{\partial}{\partial p_{-i}}$.  In this incarnation, the Heisenberg algebra is often called the \emph{Weyl algebra}.

The rich interaction between the representation theory of infinite-dimensional Lie algebras and
algebraic combinatorics owes much to the fact that the Fock space representation may be constructed in the language of symmetric functions.  Specifically, let $\Sym$ denote the algebra (over $\C$) of symmetric functions in countably many variables $\{x_1,x_2,\dots\}$.  The algebra $\Sym$ is isomorphic to a polynomial algebra in the power-sum symmetric functions:
\[
	\Sym \cong \C[P_1,P_2,\dots], \quad  P_n = \sum_{i=1}^\infty x_i^n.
\]
The algebra $\Sym$ has many vector space bases which are well-studied in the algebraic combinatorics literature, and one of the most important is the basis of Schur functions.  The natural inner product
\[
	\langle \cdot,\cdot \rangle : \Sym\times \Sym \longrightarrow \C
\]
can be defined in several ways, one of which is to declare the Schur functions to be orthonormal.  The connection with the Heisenberg algebra and the Fock space is then as follows: the Heisenberg algebra $\fh$ acts on $\Sym$, with the generator $p_{-i}$, $i \in \N_+$, acting by multiplication by the power sum $P_i$, and the generator $p_i$ of presentation~\eqref{eq:usual-h-relations'} acting by $i \frac{\partial}{\partial P_i}$, which is the linear operator adjoint to $p_{-i}$ with respect to the inner product.  This representation of $\fh$ on symmetric functions is of course isomorphic to the Fock space $\cF$; moreover, the free abelian group spanned by the Schur functions defines a natural integral form $\cF_\Z$ of the Fock space.  Any other realization of the Fock space is canonically isomorphic to $\Sym$ up to scalar multiple.  Therefore, the various natural bases of $\Sym$ often acquire interesting interpretations in other constructions of the Fock space.

The definition of the Fock space $\cF^L$ associated to a lattice Heisenberg algebra $\fh^L$ or to a quantized Heisenberg algebra is analogous to the definition of $\cF$: induce the trivial representation of the subalgebra generated by the $\{p_{i,n}\}_{n>0}$ to the entire algebra $\fh^L$.
Lattice Heisenberg algebras and their Fock space representations also have integral forms which arise naturally in geometric and categorical constructions.

%%%%%%%%%%%%%%%%%%%%%%%%%%%%%%%%%%%%%%%%%%%%%%%%%%%%%%%%%%%%%%%%%%%
%
\section{Weak categorifications and geometrizations}
\label{sec:weak-cat-and-geom}
%
%%%%%%%%%%%%%%%%%%%%%%%%%%%%%%%%%%%%%%%%%%%%%%%%%%%%%%%%%%%%%%%%%%%

The Heisenberg algebra and its Fock space representation occur as organizing objects in several places in mathematics, including the representation theory of Hecke algebras and the geometry of Hilbert schemes.  In this section we review these appearances of the Heisenberg algebra which lead to its categorification.

%%%%%%%%%%%%%%%%%%%%%%%%%%%%%%%%%%%%%%%%%%%%%%%%%%%%%%%%%%%%%%%%%%%%%%%%%%%%%%%%%%%
\subsection{The Heisenberg algebra and the representation theory of Hecke algebras}
%%%%%%%%%%%%%%%%%%%%%%%%%%%%%%%%%%%%%%%%%%%%%%%%%%%%%%%%%%%%%%%%%%%%%%%%%%%%%%%%%%%

For $n \in \N$, let $H_n(q)$ denote the Hecke algebra (over $\C$) of the symmetric group $S_n$ at a generic complex number $q$ (that is, a complex number $q$ which is not a nontrivial root of unity).  Explicitly, $H_n(q)$ is the $\C$-algebra generated by $t_1,\dots, t_{n-1}$ with relations
\begin{itemize}
  \item $t_i^2 = q + (q-1)t_i$ for $i=1,2,\dots,n$,
  \item $t_it_j = t_jt_i$ for $i,j =1,2,\dots,n-1$ such that $|j-i| > 1$,
  \item $t_it_{i+1}t_i = t_{i+1}t_it_{i+1}$ for $i=1,2,\dots,n-2$.
\end{itemize}
By convention, we set $H_0(q) = H_1(q) = \C$. To simplify notation, we will write $H_n$ for $H_n(q)$ in the sequel.

For an algebra $A$, let $A\text{-mod}$ denote the category of finite-dimensional (left) $A$-modules.  The (split) Grothendieck group of $A\text{-mod}$ is denoted $K_0(A\text{-mod})$.  This Grothendieck group is naturally a $\Z$-module, and we set
\[
	K(A\text{-mod}) = \C\otimes_\Z K_0(A\text{-mod}).
\]
The vector space $K(A\text{-mod})$ has dimension equal to the number of isomorphism classes of irreducible finite-dimensional $A$-modules.  Thus $K(H_n\text{-mod})$ has dimension equal to the number of partitions of $n$.  Following a common philosophy in the representation theory of symmetric groups, the vector spaces $K(H_n\text{-mod})$ should be studied for all $n$ at once, for when taken all together, these spaces have interesting symmetry.  Precisely, there is an isomorphism
\begin{equation}\label{eq:repiso}
	\cF \cong  \bigoplus_{n=0}^\infty K(H_n\text{-mod}) \qquad \text{(as $\fh$-modules)}
\end{equation}
between the Fock space and the direct sum of all Grothendieck groups.  This statement is essentially a theorem of Zelevinsky \cite{Zel81}, though he writes in the language of ``positive, self-adjoint Hopf algebras" and symmetric functions rather than in the language of representations of $\fh$.  When $q=1$, and $H_n$ is isomorphic to the group algebra of the symmetric group, the isomorphism~\eqref{eq:repiso} appears first in the work of Geissinger \cite{Gei77}, again in the language of the Hopf algebra of symmetric functions.

In \cite{Zel81}, Zelevinsky also considers various generalizations, wherein the Hecke algebra $H_n$ is replaced by the group algebra over $\C$ of the finite group $\GL_n(\F_q)$ or the wreath product of the symmetric group $S_n$ with an arbitrary finite group.  Zelevinsky's work was subsequently further extended by Frenkel-Jing-Wang \cite{FJW2,FJW1}, who interpreted various natural vertex operator constructions in the language of the representation theory of symmetric groups and wreath products.

We now briefly describe the action of $\fh$ on $\bigoplus_n  K(H_n\text{-mod})$ that gives rise to the isomorphism~\eqref{eq:repiso}.  The embeddings $H_n\otimes H_m \hookrightarrow H_{n+m}$ give rise to induction functors
\[
	\Ind_{n,m} :  H_n\text{-mod} \times H_m\text{-mod} \longrightarrow H_{n+m}\text{-mod}
\]
and restriction functors
\[
  \Res_{n,m} : H_{n+m}\text{-mod} \longrightarrow H_n\text{-mod} \times H_m\text{-mod}.
\]
These functors are exact and biadjoint (i.e.\ left and right adjoint) to one another.  If $M$ is a representation of $H_m$, there is an associated functor
\[
  \Ind_n(M) : H_n\text{-mod} \longrightarrow H_{n+m}\text{-mod},\quad	N \mapsto \Ind_{n,m} N\boxtimes M,
\]
which admits a biadjoint $\Res_n(M)$.

For any finite-dimensional representation $M$ of $H_m$, the functors $\Ind_n(M)$ and $\Res_n(M)$ are exact, and thus induce linear maps $K(\Ind_n(M))$ and $K(\Res_n(M))$  on the associated Grothendieck groups.  For $n \in \N_+$, define
\[
  a_n = \bigoplus_m K(\Res_m(\text{triv})), \quad b_n = \bigoplus_m K(\Ind_m(\text{sign})),
\]
where $\text{triv}$ and $\text{sign}$ are the trivial and sign representations of $H_n$.  The linear operators $a_n$ and $b_n$ satisfy the defining relations~\eqref{eq:modified-h-relations} of $\fh$ and the resulting representation of $\fh$ on $\bigoplus_n K(H_n\text{-mod})$ turns out to be isomorphic to the Fock space $\mathcal{F}$, with $\bigoplus_n K_0(H_n\text{-mod})$ isomorphic to its integral form $\mathcal{F}_\Z$.

The isomorphism~\eqref{eq:repiso} is called a \emph{weak categorification} of the Fock space because the relations in $\fh$ are shown to hold at the level of the Grothendieck group.  Strong cateorifications of $\fh$, which are the main subject of the latter parts of this survey, consist of lifting equalities in the Grothendieck group to specified isomorphisms in the category, and describing the structure of the space of natural transformations which gives rise to these isomorphisms.

%%%%%%%%%%%%%%%%%%%%%%%%%%%%%%%%%%%%%%%%%%%%%%%%%%%%%%%%%%%%%%%%%%%%%%%%%%%%%%%%%%%%%
\subsection{Heisenberg algebras and the cohomology of Hilbert schemes}
\label{subsec:coh}
%%%%%%%%%%%%%%%%%%%%%%%%%%%%%%%%%%%%%%%%%%%%%%%%%%%%%%%%%%%%%%%%%%%%%%%%%%%%%%%%%%%%%

An important motivation for the categorified Heisenberg actions to follow is the action of the Heisenberg algebra on the cohomology of Hilbert schemes.  We recall this construction briefly here, referring the reader to \cite[\S8]{Nak99} for a more complete discussion.

Let $X$ be a smooth connected (quasi-)projective surface, which, for simplicity of exposition, we assume to have cohomology only in even degree.  Let $X^{[n]}$ denote the Hilbert scheme of $n$ points on $X$.  Like symmetric groups and Hecke algebras, the Hilbert schemes $X^{[n]}$ can be studied for all $n$ together.  In particular, following G\"ottsche \cite{Go90}, in order to study the cohomology groups $\rmH^*(X^{[n]},\C)$, it is natural to consider the direct sum
$\bigoplus_{n=0}^\infty \rmH^*(X^{[n]},\C)$, which has interesting symmetry.
The intersection form, multiplied by $-1$, gives $\rmH^2(X,\Z)$ the structure of a lattice, and associated to this lattice is the Heisenberg algebra $\fh^{\rmH^2(X,\Z)}$ (see Section~\ref{sec:lattice-heisenberg}).  Constructions of Nakajima \cite{Nak99} and Grojnowski \cite{Gr96} then give an isomorphism
\begin{equation}\label{eq:geomiso}
  \cF^{\rmH^2(X,\Z)} \cong \bigoplus_{n=0}^\infty \rmH^*(X^{[n]},\C) \qquad \text{(as $\fh^{\rmH^2(X,\Z)}$-modules)}
\end{equation}
between the cohomology of Hilbert schemes and the Fock space representation $\cF^{\rmH^2(X,\Z)}$ of $\fh^{\rmH^2(X,\Z)}$.  Nakajima's construction involves defining explicit correspondences in the products of Hilbert schemes whose induced maps on cohomology satisfy the defining Heisenberg relations~\eqref{eq:usual-h-relations-lattice}.  Grojnowski's formulation is similar.  Of particular interest is the case when $X$ is the minimal resolution of a simple singularity, that is, the minimal resolution of
$\C^2/\Gamma$ for $\Gamma$ a finite subgroup of $\SL_2(\C)$, as in this case the Heisenberg algebra $\fh^{\rmH^2(X,\Z)}$ is a subalgebra of the affine Kac-Moody Lie algebra associated to $\Gamma$.  Representations of this Kac-Moody Lie algebra are often studied by first restricting to the Heisenberg algebra $\fh^{\rmH^2(X,\Z)}$, and then considering the resulting decomposition of the original representation into copies of Fock space.

The Nakajima-Grojnowski constructions are examples of \emph{geometrizations} of Heisenberg algebra representations, that is, linear actions of the Heisenberg algebra on the cohomology groups of algebraic varieties.  Geometrizations are similar to weak categorifications in that the defining relations in the algebra are only checked at the level of cohomology, rather than at the more refined level of categories of sheaves.  A principle goal of strong Heisenberg categorification is to lift the geometrizations and weak categorifications described above to strong categorifications.  In the rest of this survey we will describe several examples of such strong categorifications.

%%%%%%%%%%%%%%%%%%%%%%%%%%%%%%%%%%%%%%%%%%%%%%%%%%%%%%%%%%%%%%%%%%%
%
\section{A graphical categorification}
\label{sec:graphical-cat}
%
%%%%%%%%%%%%%%%%%%%%%%%%%%%%%%%%%%%%%%%%%%%%%%%%%%%%%%%%%%%%%%%%%%%

In this section, we review the categorification of the Heisenberg algebra in terms of planar diagrammatics introduced in \cite{LS11} (see also \cite{Kho10}) and the related strong categorification of the Fock space representation.  We choose to present these here in the language of monoidal categories instead of the 2-category theoretic language used in \cite{LS11}.

%%%%%%%%%%%%%%%%%%%%%%%%%%%%%%%%%%%
\subsection{The graphical category}
%%%%%%%%%%%%%%%%%%%%%%%%%%%%%%%%%%%

In this section, for simplicity, we will let $q$ be a complex number which is not a nontrivial root of unity, though all the constructions and theorems in this section carry over to the case where $q$ is an indeterminate.  We define an additive $\C$-linear strict monoidal category $\cH'(q)$ as follows.  The set of objects is generated by two objects $Q_+$ and $Q_-$.  Thus an arbitrary object of $\cH'(q)$ is a finite direct sum of tensor products $Q_\varepsilon := Q_{\varepsilon_1} \otimes \dots \otimes Q_{\varepsilon_n}$, where $\varepsilon = \varepsilon_1 \dots \varepsilon_n$ is a finite sequence of $+$ and $-$ signs.  The unit object $\mathbf{1}=Q_\emptyset$.

The space of morphisms $\Hom_{\cH'(q)}(Q_\varepsilon,
Q_{\varepsilon'})$ is the $\C[q,q^{-1}]$-module generated by  planar diagrams
modulo local relations.  The diagrams are
oriented compact one-manifolds immersed in the strip $\R \times
[0,1]$, modulo rel boundary isotopies.  The endpoints of the
one-manifold are located at $\{1,\dots,m\} \times \{0\}$ and
$\{1,\dots,k\} \times \{1\}$, where $m$ and $k$ are the lengths of
the sequences $\varepsilon$ and $\varepsilon'$ respectively.  The
orientation of the one-manifold at the endpoints must agree with the
signs in the sequences $\varepsilon$ and $\varepsilon'$ and triple intersections are not allowed.  For example, the diagram
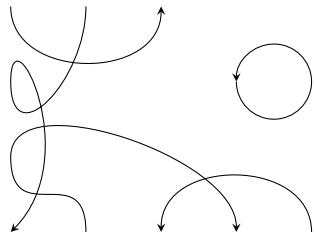
\[
  \begin{tikzpicture}[>=stealth]
    \draw[->] (0,3) .. controls (0,2) and (2,2) .. (2,3);
    \draw[->] (1,3) .. controls (1,2) and (0,1) .. (0,2) .. controls (0,3) and (1,1) .. (0,0);
    \draw[->] (1,0) .. controls (1,1) and (0,0) .. (0,1) .. controls (0,2) and (3,1) .. (3,0);
    \draw[->] (4,0) .. controls (4,1) and (2,1) .. (2,0);
    \draw[->] (3,2) arc(-180:180:.5);
  \end{tikzpicture}
\]
is a morphism from $Q_{-+--+}$ to $Q_{--+}$ (note that, in this sense, diagrams are read from bottom to top).  Composition of morphisms is given by the natural vertical gluing of diagrams and the tensor product of morphisms is horizontal juxtaposition.  An endomorphism of $\mathbf{1}$ is a diagram without endpoints.  The local relations are as follows.

\begin{equation} \label{eq:local-relation-basic-Hecke}
\begin{tikzpicture}[>=stealth,baseline=25pt]
  \draw[->] (0,0) .. controls (1,1) .. (0,2);
  \draw[->] (1,0) .. controls (0,1) .. (1,2);
  \draw (1.5,1) node {=};
  \draw (2,1) node {$q$};
  \draw[->] (2.5,0) --(2.5,2);
  \draw[->] (3.5,0) -- (3.5,2);
  \draw (4.7,1) node {$+\ (q-1)$};
  \draw[->] (5.7,0) --(6.7,2);
  \draw[->] (6.7,0) -- (5.7,2);
\end{tikzpicture}
\end{equation}

\begin{equation} \label{eq:local-relation-braid}
\begin{tikzpicture}[>=stealth,baseline=25pt]
  \draw[->] (0,0) -- (2,2);
  \draw[->] (2,0) -- (0,2);
  \draw[->] (1,0) .. controls (0,1) .. (1,2);
  \draw (2.5,1) node {=};
  \draw[->] (3,0) -- (5,2);
  \draw[->] (5,0) -- (3,2);
  \draw[->] (4,0) .. controls (5,1) .. (4,2);
\end{tikzpicture}
\end{equation}

\begin{equation} \label{eq:local-relation-up-down-double-crossing}
\begin{tikzpicture}[>=stealth,baseline=25pt]
  \draw[<-] (0,0) .. controls (1,1) .. (0,2);
  \draw[->] (1,0) .. controls (0,1) .. (1,2) ;
  \draw (1.5,1) node {=};
  \draw (2,1) node {$q$};
  \draw[<-] (2.5,0) --(2.5,2);
  \draw[->] (3.5,0) -- (3.5,2);
  \draw (4,1) node {$- \ q$};
  \draw (4.5,1.75) arc (180:360:.5) ;
  \draw (4.5,2) -- (4.5,1.75) ;
  \draw[<-] (5.5,2) -- (5.5,1.75);
  \draw (5.5,.25) arc (0:180:.5) ;
  \draw (5.5,0) -- (5.5,.25) ;
  \draw[<-] (4.5,0) -- (4.5,.25);
\end{tikzpicture} \qquad \qquad
\begin{tikzpicture}[>=stealth,baseline=25pt]
  \draw[->] (0,0) .. controls (1,1) .. (0,2);
  \draw[<-] (1,0) .. controls (0,1) .. (1,2) ;
  \draw (1.5,1) node {$=\ q$};
  \draw[->] (2.3,0) --(2.3,2);
  \draw [<-](3.3,0) -- (3.3,2);
\end{tikzpicture}
\end{equation}

\begin{equation} \label{eq:cc-circle-and-left-curl}
\begin{tikzpicture}[>=stealth,baseline=0pt]
  \draw[<-] (0,0) arc(180:0:.5);
  \draw (0,0) arc(180:360:.5);
  \draw (1.5,0) node {$=$};
  \draw (2,0) node {$\id$};
  \draw (0,-1) node {};
  \draw (0,1) node {};
\end{tikzpicture} \qquad \qquad
\begin{tikzpicture}[>=stealth,baseline=0pt]
  \draw (-1,0) .. controls (-1,.5) and (-.3,.5) .. (-.1,0) ;
  \draw (-1,0) .. controls (-1,-.5) and (-.3,-.5) .. (-.1,0) ;
  \draw (0,-1) .. controls (0,-.5) .. (-.1,0) ;
  \draw[->] (-.1,0) .. controls (0,.5) .. (0,1) ;
  \draw (0.7,0) node {$=0$};
\end{tikzpicture}
\end{equation}

\begin{defin}[Karoubi envelope] \label{def:Karoubi-envelope}
  Let $\mathcal{C}$ be a category.  The \emph{Karoubi envelope} of $\mathcal{C}$ is the category whose objects are pairs $(A,e)$ where $A$ is an object of $\mathcal{C}$ and $e \in \Hom_\mathcal{C}(A,A)$ is an idempotent endomorphism of $A$ (i.e.\ $e^2=e$).  Morphisms $(A,e) \to (A',e')$ are morphisms $f : A \to A'$ in $\mathcal{C}$ such that the diagram
  \[
    \xymatrix{
      A \ar[r]^f \ar[dr]^f \ar[d]_e & A' \ar[d]^{e'} \\
      A \ar[r]_f & A'
    }
  \]
  commutes.  Composition in the Karoubi envelope is as in $\mathcal{C}$, except that the identity morphism of $(A,e)$ is $e$.
\end{defin}

Let $\cH(q)$ be the Karoubi envelope of $\cH'(q)$.  In the case $q=1$, the category $\cH(1)$ was defined by Khovanov in \cite{Kho10}.

%%%%%%%%%%%%%%%%%%%%%%%%%%%%%%%%%%%%%%%%%%%%%%%%%%%%%%%
\subsection{Categorification of the Heisenberg algebra}
\label{subsec:graphical-h-cat}
%%%%%%%%%%%%%%%%%%%%%%%%%%%%%%%%%%%%%%%%%%%%%%%%%%%%%%%

It follows from the local relations \eqref{eq:local-relation-basic-Hecke} and \eqref{eq:local-relation-braid} that upward oriented crossings satisfy the Hecke algebra relations and so we have a canonical homomorphism
\begin{equation} \label{eq:Hn-to-Q+}
  H_n \to \End_{\cH'(q)} (Q_{+^n}).
\end{equation}
Similarly, since each space of morphisms in $\cH'(q)$ consists of diagrams up to isotopy, downward oriented crossings also satisfy the Hecke algebra relations and give us a canonical homomorphism
\begin{equation} \label{eq:Hn-to-Q-}
  H_n \to \End_{\cH'(q)} (Q_{-^n}).
\end{equation}
Introduce the complete $q$-symmetrizer and $q$-antisymmetrizer
\[ %\label{eq:q-symm-def}
  e(n) = \frac{1}{[n]_q!} \sum_{w \in S_n} t_w,\quad e'(n) = \frac{1}{[n]_{q^{-1}}!} \sum_{w \in S_n} (-q)^{-l(w)} t_w, \quad \text{where} \quad [n]_q = \sum_{i=0}^{n-1} q^i.
\]
Both $e(n)$ and $e'(n)$ are idempotents in $H_n$ (see
\cite[\S1]{Gyo86}).  We will use the notation $e(n)$ and $e'(n)$ to
also denote the image of these idempotents in $\End_{\cH'(q)}
(Q_{+^n})$ and $\End_{\cH'(q)} (Q_{-^n})$ under the canonical
homomorphisms~\eqref{eq:Hn-to-Q+} and~\eqref{eq:Hn-to-Q-}.  We then define the following objects in $\cH(q)$:
\[
  S_+^n = (Q_{+^n}, e(n)),\quad S_-^n = (Q_{-^n}, e(n)),\quad \Lambda_+^n = (Q_{+^n}, e'(n)),\quad \Lambda_-^n = (Q_{-^n}, e'(n)).
\]

\begin{theo}[{\cite[Th.~3.7]{LS11}}] \label{thm:main-object-isom}
In the category $\cH(q)$, we have
\begin{align*}
  S_-^n \otimes S_-^m &\cong S_-^m \otimes S_-^n, \\
  \Lambda_+^n \otimes \Lambda_+^m &\cong \Lambda_+^m \otimes \Lambda_+^n, \\
  S_-^n \otimes \Lambda_+^m &\cong \left( \Lambda_+^m \otimes S_-^n \right) \oplus \left( \Lambda_+^{m-1} \otimes S_-^{n-1} \right).
\end{align*}
We thus have a well-defined algebra homomorphism $\bF : \fh_\Z \to K_0(\cH(q))$ given by
\[
  \bF(a_n) = [S_-^n],\quad \bF(b_n) = [\Lambda_+^n].
\]
This homomorphism is injective.
\end{theo}

Theorem~\ref{thm:main-object-isom} was proved in \cite{Kho10} in the case $q=1$.  It is conjectured in \cite{LS11} (and in \cite{Kho10} in the case $q=1$) that $\bF$ is in fact an isomorphism.

%%%%%%%%%%%%%%%%%%%%%%%%%%%%%%%%%%%%%%%%%%%%%%%%%%%%%%%
\subsection{Categorification of Fock space}
\label{subsec:CatFock}
%%%%%%%%%%%%%%%%%%%%%%%%%%%%%%%%%%%%%%%%%%%%%%%%%%%%%%%

For $1 \le k \le n$, we can view $H_k$ as a submodule of $H_n$ via the embedding $t_i \mapsto t_i$.  We introduce here some notation for bimodules.  First note that $H_n$ is naturally an $(H_n,H_n)$-bimodule. Via our identification of $H_k$, $1 \le k \le n$, with a submodule of $H_n$, we can naturally view $H_n$ as an $(H_k,H_l)$-bimodule for $1 \le k,l \le n$.  We will write $_k (n)_l$ to denote this bimodule. If $k$ or $l$ is equal to $n$, we will often omit the subscript. Thus, for instance,
\begin{itemize}
  \item $(n)$ denotes $H_n$, considered as an $(H_n,H_n)$-bimodule,
  \item $(n)_{n-1}$ denotes $H_n$, considered as an
  $(H_n,H_{n-1})$-bimodule, and
  \item $_{n-1}(n)$ denotes $H_n$, considered as an
  $(H_{n-1},H_n)$-bimodule.
\end{itemize}
Note that tensoring on the left by $(n)_{n-1}$ (respectively $_{n-1}(n)$) is the induction functor $H_{n-1}\text{-mod} \to H_n\text{-mod}$ (respectively the restriction functor $H_n\text{-mod} \to H_{n-1}\text{-mod}$).  We have an isomorphism of $(H_n,H_n)$-bimodules
\[
  _n(n+1)_n \cong (n) \oplus \big( (n)_{n-1}(n) \big)
\]
(see \cite[Lem.~4.2]{LS11}).

We now define certain bimodule maps which will be the building blocks needed to define an action of the category $\cH(q)$ on the category of modules for Hecke algebras.  Here and in what follows, we will use string diagram notation for 2-categories.  In particular, the bimodule maps below are denoted by planar diagrams with regions labeled by elements of $\N := \Z_{\ge 0}$.  We only indicate the label of one region since the labels of the others are uniquely determined by the rule that as we move from right to left, labels increase by one as we cross upward pointing strands and decrease by one as we cross downward pointing strands.  We refer the reader to \cite{Kho10b} for an overview of this notation.

Define
\begin{equation} \label{eq:rcap-def}
  \begin{tikzpicture}[>=stealth,baseline=6pt]
    \draw (-.7,.25) node {$n+1$};
    \draw[->] (0,0) arc(180:0:.5);
  \end{tikzpicture}
  \ : (n+1)_n(n+1) \to (n+1),\quad x \otimes y \mapsto xy,
\end{equation}
to be the map given by multiplication.  Define the inclusion
\begin{equation} \label{eq:rcup-def}
  \begin{tikzpicture}[>=stealth,baseline=6pt]
    \draw (-.5,.25) node {$n$};
    \draw[->] (0,0.5) arc(180:360:.5);
  \end{tikzpicture}
  \ : (n) \hookrightarrow {_n}(n+1)_n,\quad z \mapsto z.
\end{equation}
Define
\begin{equation} \label{eq:lcap-def}
  \begin{tikzpicture}[>=stealth,baseline=6pt]
    \draw (-.5,.25) node {$n$};
    \draw[<-] (0,0) arc(180:0:.5);
  \end{tikzpicture}
  \ : {_n}(n+1)_n \to (n)
\end{equation}
to be the map that is the identity on $H_n$ and that maps $t_n$ to zero.  Finally, define
\begin{equation} \label{eq:lcup-def}
  \begin{tikzpicture}[>=stealth,baseline=6pt]
    \draw (-.7,.25) node {$n+1$};
    \draw[<-] (0,0.5) arc(180:360:.5);
  \end{tikzpicture}
  \ : (n+1) \to (n+1)_n(n+1)
\end{equation}
to be the map determined by
\[
  1_{n+1} \mapsto \sum_{i=1}^{n+1} q^{i-(n+1)} t_{i}\hdots t_{n-1} t_n \otimes t_n t_{n-1}\hdots t_{i}.
\]
We set $t_{n+1}=1$ above, so that the $i=n+1$ term in the sum is $1\otimes 1$.  The maps defined in \eqref{eq:rcap-def}--\eqref{eq:lcup-def} are adjunction maps that make $(\Res,\Ind)$ into a biadjoint pair (see \cite[Prop.~4.4]{LS11}).

Our final diagrammatic ingredient is the crossing, which we define as follows:
\begin{equation} \label{eq:upcross-def}
  \begin{tikzpicture}[>=stealth,baseline=10pt]
    \draw[->] (0,0) to (1,1);
    \draw[->] (1,0) to (0,1);
    \draw (1.1,.5) node {$n$};
  \end{tikzpicture}
  \quad : \quad (n+2)_n \to (n+2)_n,\quad z \mapsto
  z t_{n+1}.
\end{equation}
It follows from \cite[Prop.~4.5]{LS11} that any two isotopic diagrams involving this crossing as well as cups and caps will give rise to the same bimodule map (see below).

For $n,m \in \N$, let $_m \bimod_n$ be the category of finite-dimensional $(H_m,H_n)$-bimodules and let $\bimod_n = \bigoplus_{m \in \N} {_m \bimod_n}$.  By \cite[Prop.~4.7]{LS11}, the relations \eqref{eq:local-relation-basic-Hecke}--\eqref{eq:cc-circle-and-left-curl} hold when these diagrams are interpreted as maps of bimodules (with any labelings of the regions).  For $n \in \N$, we therefore have a well-defined functor $\bA_n : \cH'(q) \to \bimod_n$ defined as follows.  For an object $Q_\varepsilon$ of $\cH'(q)$, $\varepsilon = \varepsilon_1 \varepsilon_2 \dots \varepsilon_\ell$, $\bA(Q_\varepsilon)$ is the tensor product of induction of restriction bimodules, where each $+$ corresponds to the induction bimodule and each $-$ corresponds to the restriction bimodule.  For instance,
\[
  \bA_n(Q_{+--+-++}) = (n+1)_{n} (n+1)_{n+1} (n+2)_{n+2} (n+2)_{n+1} (n+2)_{n+2} (n+2)_{n+1} (n+1)_n.
\]
In remains to define $\bA_n$ on planar diagrams (morphisms of $\cH'(q)$).  This is done as follows.  Let $D$ be a morphism of $\cH'(q)$.  We label the rightmost region of the diagram $D$ with $n$.  We then label all other regions of the diagram with integers such that as we move across the diagram from right to left, labels increase by one when we cross upward pointing strands and decrease by one when we cross downward pointing strands.  It is easy to see that there is a unique way to do this.  For instance, the following diagram would be labeled as indicated.
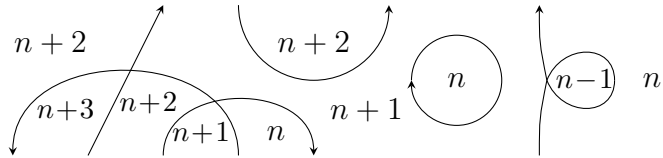
\[
  \begin{tikzpicture}[>=stealth]
    \draw[->] (0,0) to (1,2);
    \draw[->] (1,0) .. controls (1,1) and (3,1) .. (3,0);
    \draw[->] (2,0) .. controls (2,1.5) and (-1,1.5) .. (-1,0);
    \draw[->] (2,2) arc(180:360:1);
    \draw (7,1) .. controls (7,1.5) and (6.3,1.5) .. (6.1,1);
    \draw (7,1) .. controls (7,.5) and (6.3,.5) .. (6.1,1);
    \draw (6,0) .. controls (6,.5) .. (6.1,1);
    \draw[->] (6.1,1) .. controls (6,1.5) .. (6,2);
    \draw[->] (4.3,1) arc(180:-180:.6);
    \draw (7.5,1) node {$n$};
    \draw (6.6,1) node {\small{$n\!-\!1$}};
    \draw (3.7,.6) node {$n+1$};
    \draw (4.9,1) node {$n$};
    \draw (3,1.5) node {$n+2$};
    \draw (-.5,1.5) node {$n+2$};
    \draw (-.3,.6) node {\small$n\!+\!3$};
    \draw (.8,.7) node {\small$n\!+\!2$};
    \draw (1.5,.3) node {\small$n\!+\!1$};
    \draw (2.5,.3) node {$n$};
  \end{tikzpicture}
\]
The functor $\bA_n$ then maps the labeled diagram $D$ to the corresponding bimodule map according to the definitions~\eqref{eq:rcap-def}--\eqref{eq:upcross-def}.  More precisely, we isotope $D$ so that it contains left and right cups and caps in addition to upward pointing crossings (alternatively, one can directly compute the bimodule maps corresponding to the other crossings and then there is no need to isotope $D$).  Then $D$ is mapped to the corresponding composition of the maps~\eqref{eq:rcap-def}--\eqref{eq:upcross-def}.  We adopt the convention that $D$ is mapped to zero if any of its regions is labeled by a negative integer.  Since the category $\bimod_n$ is idempotent complete, the functor $\bA_n$ induces a functor $\bA_n : \cH(q) \to \bimod_n$ which we denote by the same symbol.

For abelian or triangulated categories $\mathcal{C}, \mathcal{D}$, let $\Fun(\mathcal{C},\mathcal{D})$ denote the category of exact functors from $\mathcal{C}$ to $\mathcal{D}$ (with morphisms being natural transformations).  Any element of $\Fun(\mathcal{C},\mathcal{D})$ induces a $\Z$-linear map $K_0(\mathcal{C}) \to K_0(\mathcal{D})$.  We thus have a natural map
\[
  K_0 : \Fun(\mathcal{C},\mathcal{D}) \to \Hom_\Z(K_0(\mathcal{C}), K_0(\mathcal{D})).
\]

Now, for $m,n \in \N$, there is natural functor from $_m \bimod_n$ to $\Fun(H_n\text{-mod},H_m\text{-mod})$.  Namely, an $(H_m,H_n)$-bimodule $M$ is sent to the functor $M \otimes_{H_n} \cdot$ (i.e.\ the functor given by tensoring on the left with $M$) and bimodule maps are sent to the corresponding natural transformations.  We then define the functor $\bA$ to be the composition
\begin{multline} \label{eq:functor-A-def}
  \bA : \cH(q) \xrightarrow{\bigoplus_{n \in \N} \bA_n} \bigoplus_{m,n \in \N} {_m \bimod_n} \to \bigoplus_{m,n \in \N} \Fun (H_n\text{-mod},H_m\text{-mod}) \\
  \xrightarrow{} \Fun \left( \bigoplus_{n \in \N} H_n\text{-mod}, \bigoplus_{n \in \N} H_n\text{-mod} \right).
\end{multline}
Thus $\bA$ defines a representation of the category $\cH(q)$ on the category $\bigoplus_{n \in \N} H_n\text{-mod}$.

\begin{theo}[Categorification of Fock space]\label{thm:Fock}
We have a commutative diagram
\begin{equation} \label{eq:Fock-space-cat-diagram}
  \vcenter{\xymatrix{
    \cH(q) \ar[r]^(0.18){\bA} \ar[d]_{K_0} & \Fun \left( \bigoplus_{n \in \N} H_n\textup{-mod}, \bigoplus_{n \in \N} H_n\textup{-mod} \right) \ar[d]^{K_0} \\
    K_0(\cH(q)) \ar[r]^{K_0(\bA)} & \End_\Z \mathcal{F}_\Z \\
    \fh_\Z \ar@{^{(}->}[u]^{\bF} \ar[ur]
  }},
\end{equation}
where the arrow $\fh_\Z \to \End_\Z \mathcal{F}_\Z$ is the action of the integral version of the Heisenberg algebra $\fh_\Z$ on the integral version $\mathcal{F}_\Z$ of the Fock space representation described in Section~\ref{sec:weak-cat-and-geom}.  Thus the action of $\cH(q)$ on $\bigoplus_{n \in \N} H_n\textup{-mod}$ is a (strong) categorification of the Fock space representation of the Heisenberg algebra (one can recover the usual, non-integral version, after tensoring with $\C$ or by replacing $K_0$ by $K$ everywhere in the above).
\end{theo}

There is an isomorphism from $\mathcal{F}_\Z \cong \Z[b_1,b_2,\dots]$ to the ring of symmetric functions, taking $b_i$ to the elementary symmetric function, usually denoted $e_i$.  Theorem \ref{thm:Fock} endows the Fock space with another natural basis, namely, the classes of irreducible $H_n$-modules in the Grothendieck group.  One can show that under the isomorphism to symmetric functions, these classes get mapped to the basis of Schur functions.

As is usual in the theory of categorification, the category $\cH(q)$ has considerably more structure than the ring $\fh_\Z$ which it categorifies.  Much of this extra structure appears in the Hom-spaces of $\cH(q)$.  It turns out that one can give an explicit description of these Hom-spaces.

In order to simplify our pictures, we will use a hollow dot to denote a right curl and a hollow dot labeled $d$ to denote $d$ right curls.
\begin{equation} \label{eq:right-curl-shorthand}
\begin{tikzpicture}[>=stealth,baseline=25pt]
  \draw[->] (0,0) to (0,2);
  \draw [red] (0,1) circle (4pt);
  \draw (0.5,1) node {$:=$};
  \draw (2,1) .. controls (2,1.5) and (1.3,1.5) .. (1.1,1);
  \draw (2,1) .. controls (2,.5) and (1.3,.5) .. (1.1,1);
  \draw (1,0) .. controls (1,.5) .. (1.1,1);
  \draw[->] (1.1,1) .. controls (1,1.5) .. (1,2);
\end{tikzpicture}
\qquad \qquad
\begin{tikzpicture}[>=stealth,baseline=25pt]
  \draw[->] (0,0) to (0,2);
  \draw[red] (0,1) circle (4pt);
  \draw (0.1,1) node[anchor=west] {$d$};
  \draw (.85,1) node {$:=$};
  \draw[->] (1.5,0) to (1.5,2);
  \draw[red] (1.5,.45) circle (4pt);
  \draw[red] (1.5,.8) circle (4pt);
  \draw[red] (1.5,1.15) circle (4pt);
  \draw[red] (1.5,1.5) circle (4pt);
  \draw (2.5,1) node {$\Bigg\}\ d$ dots};
\end{tikzpicture}
\end{equation}
We note that in \cite{Kho10,LS11}, right curls were denoted by solid dots.  We choose to use hollow dots here to match the notation of \cite{CauLic10} and used in Section~\ref{sec:quantum-heis-cat}.  Define $c_d$ to be a clockwise oriented circle with $d$ right curls.
\[
\begin{tikzpicture}[>=stealth]
  \draw (-2,0) node {$c_d:=$};
  \draw[<-] (0,0) arc (0:360:.5);
  \draw[red] (-1,0) circle (4pt);
  \draw (-1.05,0) node[anchor=east] {$d$};
\end{tikzpicture}
\]

Let $H_n^+$ be the $\C$-algebra with generators $t_i, y_i$, $1 \le i \le n$, and defining relations
\begin{align*}
  y_it_k &= t_ky_i,\quad i\neq k, k+1, \\
  t_i y_{i+1} &= y_it_i + (q-1)y_{i+1} + q, \\
  y_{i+1}t_i &= t_iy_i + (q-1)y_{i+1} + q.
\end{align*}
If $q \neq 1$, $H_n^+$ is isomorphic to a natural subalgebra of the affine Hecke algebra of type $A$ (see \cite[Lem.~2.7]{LS11}).  Moreover, when $q=1$, $H_n^+$ is isomorphic to the \emph{degenerate affine Hecke algebra} of type $A$.

It follows from \cite[Lem.~2.3]{LS11} that there is a natural morphism
\[
  \phi_n:H_n^+\to \End_{\cH(q)}(Q_{+^n})
\]
taking $t_k$ to the crossing of the $k$ and $(k+1)$-st strands and taking $y_k$ to a right curl (or hollow dot) on the $k$-th strand.  More generally, there is a natural morphism
\[
  \psi_n = \phi_n\otimes \psi_0 : H_n^+ \otimes_{\C[q,q^{-1}]} \C[q,q^{-1}][c_0,c_1,\dots]
  \to \End_{\cH'(q)}(Q_{+^n}),
\]
where the dotted clockwise circles corresponding to elements of $\C[c_0,c_1,\dots]$ are placed to the right of the diagrams corresponding to elements of $H_n^+$.  By \cite[Th.~2.8]{LS11}, $\psi_n$ is in fact an isomorphism of algebras, thus giving a precise description of the space $\End_{\cH'(q)}(Q_{+^n})$.

It is also possible to give an explicit basis for an arbitrary Hom-space $\Hom_{\cH'(q)}(Q_\varepsilon, Q_{\varepsilon'})$ for any sequences $\varepsilon, \varepsilon'$.  Let $k$ denote the total number of $+$s in $\epsilon$ and $-$s in $\epsilon'$.  We clearly have $\Hom_{\cH'(q)}(Q_\epsilon, Q_{\epsilon'})=0$ if the total number of $-$s in $\epsilon$ and $+$s in $\epsilon'$ is not also equal to $k$.  Thus, we assume from now on that $k$ is also the total number of $-$s in $\epsilon$ and $+$s in $\epsilon'$.

\begin{prop}[{\cite[Prop.~2.10]{LS11}}]
For any sign sequences $\varepsilon, \varepsilon'$, a basis of the $\C[q,q^{-1}]$-module $\Hom_{\cH'(q)} (Q_\varepsilon,Q_{\varepsilon'})$ is given by the set $B(\varepsilon,\varepsilon')$, which is the set of planar diagrams obtained in the following manner:
\begin{itemize}
  \item The sequences $\varepsilon$ and $\varepsilon'$ are written at the bottom and top (respectively) of the plane strip $\R \times [0,1]$.

  \item The elements of $\varepsilon$ and $\varepsilon'$ are matched by oriented segments embedded in the strip in such a way that their orientations match the signs (that is, they start at either a $+$ of $\epsilon$ or a $-$ of $\epsilon'$, and end at either a $-$ of $\epsilon$ or a $+$ of $\epsilon'$), each two segments intersect at most once, and no triple intersections are allowed.

  \item Any number of hollow dots may be placed on each interval near its out endpoint (i.e. between its out endpoint and any intersections with other intervals).

  \item In the rightmost region of the diagram, a finite number of clockwise disjoint nonnested circles with any number of hollow dots may be drawn.
\end{itemize}
Thus the set of diagrams $B(\varepsilon,\varepsilon')$ is parameterized by $k!$ possible matchings of the $2k$ oriented endpoints, a sequence of $k$ nonnegative integers determining the number of hollow dots on each interval, and by a finite sequence of nonnegative integers determining the number of clockwise circles with various numbers of hollow dots.
\end{prop}

An example of an element of $B(--+-+,+-+-+--)$ is drawn below.
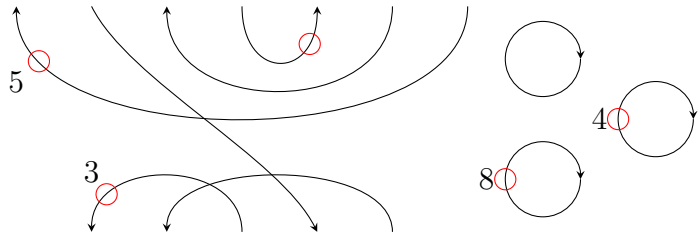
\[
\begin{tikzpicture}[>=stealth]
  \draw[<-] (0,0) .. controls (0,1) and (2,1) .. (2,0);
  \draw[<-] (1,0) .. controls (1,1) and (4,1) .. (4,0);
  \draw[->] (0,3) .. controls (.5,2) and (2.5,1) .. (3,0);
  \draw[<-] (-1,3) .. controls (-1,1) and (5,1) .. (5,3);
  \draw[<-] (1,3) .. controls (1,1.5) and (4,1.5) .. (4,3);
  \draw[->] (2,3) .. controls (2,2) and (3,2) .. (3,3);
  \draw[red] (-.7,2.27) circle (4pt);
  \draw (-1,2) node {$5$};
  \draw[red] (.2,.5) circle (4pt);
  \draw (0,.8) node {$3$};
  \draw[red] (2.9,2.5) circle (4pt);
  \draw[<-] (6.5,2.3) arc(0:360:.5);
  \draw[<-] (6.5,.7) arc(0:360:.5);
  \draw[red] (5.5,.7) circle (4pt);
  \draw (5.5,.7) node[anchor=east] {$8$};
  \draw[<-] (8,1.5) arc(0:360:.5);
  \draw[red] (7,1.5) circle (4pt);
  \draw (7,1.5) node[anchor=east] {$4$};
\end{tikzpicture}
\]

It turns out that the hollow dots (right curls) have a nice interpretation in terms of the Fock space categorification via modules for Hecke algebras.  For $k = 0,1,2,\hdots n$, let
\begin{align*}
  L_{k+1} &= \sum_{i=1}^k q^{i-k} t_i \cdots t_k \cdots t_i \\
  &= t_{k} + q^{-1} t_{k-1}t_{k}t_{k-1} +  q^{-2} t_{k-2}t_{k-1}t_kt_{k-1}t_{k-2}
  + \dots + q^{1-k} t_1 \cdots t_k \cdots t_1.
\end{align*}
By convention, $L_1=0$.  The $L_k$ (or, more precisely, $q^{-1}L_k$) are called \emph{Jucys-Murphy elements} of $H_{n+1}$ (see, for example, \cite[\S3.3]{Mat99}).  By \cite[Prop.~4.12]{LS11}, the functor $\bA$ (see~\eqref{eq:functor-A-def}) maps the single right curl
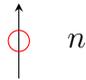
\[
\begin{tikzpicture}[>=stealth]
  \draw[->] (0,0) to (0,1);
  \draw[red] (0,.5) circle (4pt);
  \draw (0.5,.5) node[anchor=west] {$n$};
\end{tikzpicture}
\]
to the bimodule map
\[
  (n+1)_n \to (n+1)_n,\quad z \mapsto zL_{n+1}.
\]
Thus right curls correspond to Jucys-Murphy elements, which play a significant role in the theory of Hecke algebras.  The observation that right curls correspond to Jucys-Murphy elements was first made by Khovanov in the case $q=1$ (see~\cite[\S4]{Kho10}).

The category $\cH(q)$ also acts on the category $\bigoplus_n \C[\GL_n(\mathbb{F}_q)]$-mod of modules over finite general linear groups.  In this case, upward and downward oriented strands correspond to certain parabolic induction and restriction functors.  This action yields a second categorification of the Fock space representation of $\fh_\Z$.  We refer the reader to \cite[\S5]{LS11} for details.

%%%%%%%%%%%%%%%%%%%%%%%%%%%%%%%%%%%%%%%%%%%%%%%%%%%%%%%%%%%%%%%%%%%
%
\section{Quantum Heisenberg categorification, finite subgroups of $\SL_2(\C)$, and Hilbert schemes}
\label{sec:quantum-heis-cat}
%
%%%%%%%%%%%%%%%%%%%%%%%%%%%%%%%%%%%%%%%%%%%%%%%%%%%%%%%%%%%%%%%%%%%

In this section we describe the strong categorification, developed in \cite{CauLic10}, that is related to the geometrization via the Hilbert scheme described in Section~\ref{subsec:coh}.

Let $\Gamma$ be a finite subgroup of $\SL_2(\C)$.  To each conjugacy class of finite subgroup one can associate an affine Dynkin diagram, which we also denote by $\Gamma$, and hence one of the lattices $L_\Gamma$ from Section~\ref{sec:heisenberg-defs}. This correspondence is due originally to McKay \cite{Mc80}. We will review a slightly modified version of this correspondence, which constructs a $t$-deformation of the lattice $L_\Gamma$ from the finite group $\Gamma$.

The inclusion $\Gamma \subseteq \SL_2(\C)$ defines a 2-dimensional representation $V$ of $\Gamma$. As a result, $\Gamma$ acts by automorphisms on the exterior algebra $\bigwedge^*(V)$,  and one can form the smash product algebra $B^\Gamma := \C[\Gamma]\#\bigwedge^*(V)$.  The algebra
$B^\Gamma$ is a finite-dimensional Frobenius algebra with a nondegenerate trace
$\tr:B^\G \to \C$.  This algebra inherits a $\Z$-grading from the natural grading on the exterior algebra.

Let $B^\Gamma\text{-gpmod}$ be the category of finite-dimensional graded projective $B_\Gamma$-modules, and let $K_0(B_\Gamma\text{-gpmod})$ be its (split) Grothendieck group.  Because we are considering graded modules, $K_0(B_\Gamma\text{-gpmod})$ is actually a $\Z[t,t^{-1}]$ module, where multiplication by $t$ and $t^{-1}$ on the Grothendieck group come from positive and negative grading shifts in the category $B^\Gamma\text{-gpmod}$,
\[
	t^{\pm 1}[M] := [M\langle \pm 1 \rangle].
\]
The $\Hom$ bifunctor
\[
	\Hom : B^\Gamma\text{-gpmod} \times B^\Gamma\text{-gpmod} \longrightarrow \text{gVect},
\]
whose image is the category $\text{gVect}$ of finite-dimensional $\Z$-graded vector spaces,
descends to a $\Z[t,t^{-1}]$ semilinear pairing on $K_0(B_\Gamma\text{-gpmod})$:
\[
	\langle [M],[N]\rangle = \dim_t(\Hom_{B_\Gamma} (M,N)) \in K_0(\text{gVect}) \cong \Z[t,t^{-1}].
\]
In the above, if $W = \bigoplus_{n \in \Z} W(n)$ is a $\Z$-graded vector space, then $\dim_t(W)$ is its graded dimension:
\[
	\dim_t(W) = \sum_{n \in \Z} t^n \dim_{\C} W(n).
\]
This pairing is a $t$-deformation of the bilinear form of Section~\ref{sec:heisenberg-defs}.

%%%%%%%%%%%%%%%%%%%%%%%%%%%%%%%%%%%%%%%%%%%%%%%%%%%%%%%%%%%%%%%%%%%%%%%%%%%%%%%%
\subsection{Categorification of the quantized Heisenberg algebra}
%%%%%%%%%%%%%%%%%%%%%%%%%%%%%%%%%%%%%%%%%%%%%%%%%%%%%%%%%%%%%%%%%%%%%%%%%%%%%%%%

Recall the presentation of the quantized Heisenberg algebra $\fh^{L_{\G,t}}$ associated to the affine Dynkin diagram $\Gamma$ from Section \ref{sec:lattice-heisenberg}.  In this section, we recall the graphical calculus giving rise to a strong categorification of $\fh^{L_{\G,t}}$ constructed in \cite{CauLic10}.

We define an additive strict monoidal category $\cH'_\G$ as follows.  The set of objects of $\cH'_\G$ is generated by symbols $Q_+\langle n\rangle $ and $Q_-\langle n\rangle$, for $n\in \Z$ (the symbols  $Q_+$ and $Q_-$ were denoted by $P$ and $Q$, respectively, in \cite{CauLic10}).  Just as in Section~\ref{sec:graphical-cat}, the monoidal
unit object is denoted $\mathbf{1}=Q_\emptyset \langle 0 \rangle$.  Moreover, the modoidal structure is declared to be compatible with shifts $\langle \cdot \rangle$ in the sense that $Q_\varepsilon \langle s \rangle \otimes Q_{\varepsilon'} \langle s' \rangle = Q_{\varepsilon \varepsilon'} \langle s+s' \rangle$.  Thus an arbitrary object of $\cH'_\G$ is a finite direct sum of tensor products $Q_\varepsilon \langle s \rangle := Q_{\varepsilon_1}\langle s_1\rangle \otimes
\dots \otimes Q_{\varepsilon_n}\langle s_n \rangle$, where $\varepsilon = \varepsilon_1
\dots \varepsilon_n$ is a finite sequence of $+$ and $-$ signs, and $s_1+ \dots +s_n = s$.

The space of morphisms $\Hom_{\cH'_\G}(Q_\varepsilon\langle s \rangle, Q_{\varepsilon'}\langle s'\rangle)$ is the vector space generated by planar diagrams modulo local relations.  The diagrams are oriented compact one-manifolds immersed in the strip $\R \times [0,1]$, modulo rel boundary isotopies.  Each such diagram has a grading, to be defined later in this section, which determines the difference in shift in the domain and codomain.

In a given planar diagram, the endpoints of the one-manifold are located at $\{1,\dots,m\} \times \{0\}$ and $\{1,\dots,k\} \times \{1\}$, where $m$ and $k$ are the lengths of the sequences $\varepsilon$ and $\varepsilon'$ respectively (this is very similar to the setup of Section~\ref{sec:graphical-cat}).  The orientation of the one-manifold at the endpoints must agree with the signs in the sequences $\varepsilon$ and $\varepsilon'$, and triple intersections are not allowed.  A new feature appearing here which did not appear in Section~\ref{sec:graphical-cat} is that each immersed one-manifold is allowed to carry dots labeled by elements $b \in B^\G$.  For example,
\[
  \begin{tikzpicture}[>=stealth]
  \draw[->] (0,0) -- (0,1.5) [];
  \filldraw [blue](0,.375) circle (2pt);
  \draw (0,.375) node [anchor=west] [black] {$b''$};
  \filldraw [blue](0,.75) circle (2pt);
  \draw (0,.75) node [anchor=west] [black] {$b'$};
  \filldraw [blue](0,1.125) circle (2pt);
  \draw (0,1.125) node [anchor=west] [black] {$b$};
  \end{tikzpicture}
\]
is an element of $\Hom_{\cH'_\G}(Q_+\langle s\rangle ,Q_+\langle s+ \mbox{deg}(b)+\mbox{deg}(b') + \mbox{deg}(b'')\rangle)$ for every $s\in \Z$, while
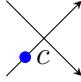
\[
  \begin{tikzpicture}[>=stealth]
  \draw[->] (0,0) -- (1,1) [];
  \draw[<-] (1,0) -- (0,1) [];
  \filldraw [blue](.25,.25) circle (2pt);
  \draw (.25,.25) node [anchor=west] [black] {$c$};
  \end{tikzpicture}
\]
is an element of $\Hom_{\cH'_\G}(Q_{+-}\langle s \rangle,Q_{-+}\langle s + \mbox{deg}(c)\rangle)$ for every $s\in \Z$. Note that the domain of a morphism is specified (up to shift) at the bottom of the diagram, the codomain is specified (up to shift) at the top, and compositions of morphisms are read from bottom to top.   That the difference between the shifts in the domain and codomain is $ \mbox{deg}(b)+\mbox{deg}(b') + \mbox{deg}(b'')$ in the first picture and $ \mbox{deg}(c)$ in the second picture will be clear once the degrees of these diagrams have been defined.

The local relations imposed are the following. First we have relations involving the movement of dots along the carrier strand. We allow dots to move freely along strands and through intersections:
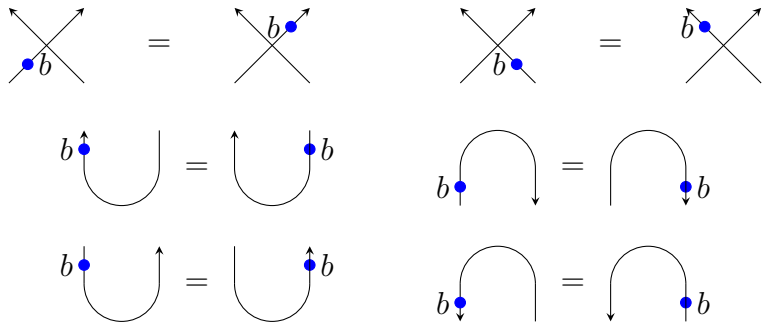
\[
  \begin{tikzpicture}[>=stealth]
  \draw[->] (0,0) -- (1,1) [];
  \draw[->] (1,0) -- (0,1) [];
  \filldraw [blue](.25,.25) circle (2pt);
  \draw (.25,.25) node [anchor=west] [black] {$b$};
  \draw (2,.5) node{$= $};
  \draw [shift={+(1,0)}][->](2,0) -- (3,1) [];
  \draw [shift={+(1,0)}][->](3,0) -- (2,1) [];
  \filldraw [shift={+(1,0)}][blue](2.75,0.75) circle (2pt);
  \draw [shift={+(1,0)}](2.75,.75) node [anchor=east] [black] {$b$};
  \draw  [shift={+(1,0)}][shift={+(5,0)}][->](0,0) -- (1,1) [];
  \draw  [shift={+(1,0)}][shift={+(5,0)}][->](1,0) -- (0,1) [];
  \filldraw  [shift={+(1,0)}][shift={+(5,0)}][blue](.75,.25) circle (2pt);
  \draw [shift={+(1,0)}][shift={+(5,0)}] (.75,.25) node [anchor=east] [black] {$b$};
  \draw  [shift={+(1,0)}][shift={+(5,0)}](2,.5) node{$= $};
  \draw[->]  [shift={+(2,0)}][shift={+(5,0)}](2,0) -- (3,1) [];
  \draw[->]  [shift={+(2,0)}][shift={+(5,0)}](3,0) -- (2,1) [];
  \filldraw  [shift={+(2,0)}][shift={+(5,0)}][blue](2.25,0.75) circle (2pt);
  \draw [shift={+(2,0)}][shift={+(5,0)}] (2.25,.75) node [anchor=east] [black] {$b$};
  \end{tikzpicture}
\]
\[
  \begin{tikzpicture}[>=stealth]
  \draw[->] (0,0) -- (0,.5);
  \filldraw [blue] (0,.25) circle (2pt);
  \draw (0,.25) node [anchor=east] [black] {$b$};
  \draw (1,0) -- (1,.5);
  \draw (0,0) arc (180:360:.5);
  \draw (1.5,0) node {=};
  \draw[->] (2,0) -- (2,.5);
  \filldraw [blue] (3,.25) circle (2pt);
  \draw (3,.25) node [anchor=west] [black] {$b$};
  \draw (3,0) -- (3,.5);
  \draw (2,0) arc (180:360:.5);

  \draw (5,0) -- (5,-.5);
  \filldraw [blue] (5,-.25) circle (2pt);
  \draw (5,-.25) node [anchor=east] [black] {$b$};
  \draw[->] (6,0) -- (6,-.5);
  \draw (5,0) arc (180:0:.5);
  \draw (6.5,0) node {=};
  \draw (7,0) -- (7,-.5);
  \filldraw [blue] (8,-.25) circle (2pt);
  \draw (8,-.25) node [anchor=west] [black] {$b$};
  \draw[->] (8,0) -- (8,-.5);
  \draw (7,0) arc (180:0:.5);
\end{tikzpicture}
\]
\[
\begin{tikzpicture}[>=stealth]
\draw (0,0) -- (0,.5);
\filldraw [blue] (0,.25) circle (2pt);
\draw (0,.25) node [anchor=east] [black] {$b$};
\draw[->] (1,0) -- (1,.5);
\draw (0,0) arc (180:360:.5);
\draw (1.5,0) node {=};
\draw (2,0) -- (2,.5);
\filldraw [blue] (3,.25) circle (2pt);
\draw (3,.25) node [anchor=west] [black] {$b$};
\draw[->] (3,0) -- (3,.5);
\draw (2,0) arc (180:360:.5);

\draw[->] (5,0) -- (5,-.5);
\filldraw [blue] (5,-.25) circle (2pt);
\draw (5,-.25) node [anchor=east] [black] {$b$};
\draw (6,0) -- (6,-.5);
\draw (5,0) arc (180:0:.5);
\draw (6.5,0) node {=};
\draw[->] (7,0) -- (7,-.5);
\filldraw [blue] (8,-.25) circle (2pt);
\draw (8,-.25) node [anchor=west] [black] {$b$};
\draw (8,0) -- (8,-.5);
\draw (7,0) arc (180:0:.5);
\end{tikzpicture}
\]
The U-turn 2-morphisms (i.e.\ the left and right cups and caps) are adjunctions making $Q_+$ and $Q_-$ biadjoint up to a grading shift. Collision of dots is controlled by multiplication in the algebra $B^\G$:
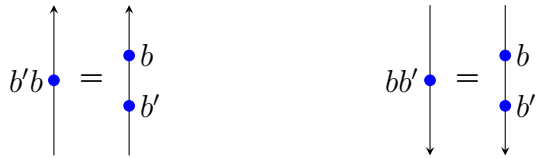
\[
\begin{tikzpicture}[>=stealth]
\draw[->] (0,0) -- (0,2);
\filldraw [blue] (0,1) circle (2pt);
\draw (0,1) node [anchor=east] [black] {$b'b$};
\draw (.5,1) node {=};
\draw (.5,1) node {=};
\draw[->] (1,0) -- (1,2);
\filldraw [blue] (1,.66) circle (2pt);
\draw (1,.66) node [anchor=west] [black] {$b'$};
\filldraw [blue] (1,1.33) circle (2pt);
\draw (1,1.33) node [anchor=west] [black] {$b$};

\draw[->] (5,2) -- (5,0);
\filldraw [blue] (5,1) circle (2pt);
\draw (5,1) node [anchor=east] [black] {$bb'$};
\draw (5.5,1) node {=};
\draw (5.5,1) node {=};
\draw[->] (6,2) -- (6,0);
\filldraw [blue] (6,.66) circle (2pt);
\draw (6,.66) node [anchor=west] [black] {$b'$};
\filldraw [blue] (6,1.33) circle (2pt);
\draw (6,1.33) node [anchor=west] [black] {$b$};
\end{tikzpicture}
\]
Dots on distinct strands supercommute when they move past one another:
\[
  \begin{tikzpicture}[>=stealth]
  \draw[->] (0,0) -- (0,2);
  \filldraw [blue] (0,.66) circle (2pt);
  \draw (0,.66) node [anchor=east] [black] {$b$};
  \draw (.5,.5) node {$\hdots$};
  \draw[->] (1,0)--(1,2);
  \filldraw [blue] (1,1.33) circle (2pt);
  \draw (1,1.33) node [anchor=west] [black] {$b'$};
  \draw (3.5,1) node {$= \ (-1)^{(\deg b)\cdot(\deg b')}$};
  \draw[->] [shift={+(3.7,0)}](2,0) --(2,2);
  \filldraw [shift={+(3.7,0)}][blue] (2,1.33) circle (2pt);
  \draw [shift={+(3.7,0)}](2,1.33) node [anchor=east] [black] {$b$};
  \draw [shift={+(3.7,0)}](2.5,.5) node {$\hdots$};
  \filldraw [shift={+(3.7,0)}][blue] (3,.66) circle (2pt);
  \draw [shift={+(3.7,0)}](3,.66) node [anchor=west] [black] {$b'$};
  \draw[->] [shift={+(3.7,0)}](3,0) -- (3,2);
  \end{tikzpicture}
\]
In addition to specifying how dots collide and slide we also impose the following local relations:
\begin{equation}\label{eq:rel1}
  \begin{tikzpicture}[>=stealth,baseline=25pt]
  \draw[->] [shift={+(7,0)}](0,0) .. controls (1,1) .. (0,2);
  \draw[->] [shift={+(7,0)}](1,0) .. controls (0,1) .. (1,2);
  \draw [shift={+(7,0)}](1.5,1) node {=};
  \draw[->] [shift={+(7,0)}](2,0) --(2,2);
  \draw[->] [shift={+(7,0)}](3,0) -- (3,2);

  \draw[->] (0,0) -- (2,2);
  \draw[->] (2,0) -- (0,2);
  \draw[->] (1,0) .. controls (0,1) .. (1,2);
  \draw (2.5,1) node {=};
  \draw[->] (3,0) -- (5,2);
  \draw[->] (5,0) -- (3,2);
  \draw[->] (4,0) .. controls (5,1) .. (4,2);
  \end{tikzpicture}
\end{equation}

\begin{equation}\label{eq:rel2}
  \begin{tikzpicture}[>=stealth,baseline=25pt]
  \draw[<-][shift={+(-1,0)}] (0,0) .. controls (1,1) .. (0,2);
  \draw[->][shift={+(-1,0)}] (1,0) .. controls (0,1) .. (1,2);
  \draw[shift={+(-1,0)}] (1.5,1) node {=};
  \draw[<-][shift={+(-1,0)}] (2,0) --(2,2);
  \draw[->][shift={+(-1,0)}] (3,0) -- (3,2);

  \draw[shift={+(-0.6,0)}] (3.8,1) node{$-\sum_{b \in \mathcal{B}}$};

  \draw (4,1.75) arc (180:360:.5);
  \draw (4,2) -- (4,1.75);
  \draw[<-] (5,2) -- (5,1.75);
  \draw (5,.25) arc (0:180:.5);
  \filldraw [blue] (4.5,1.25) circle (2pt);
  \draw (4.5,1.25) node [anchor=south] {$b$};
  \filldraw [blue] (4.5,0.75) circle (2pt);
  \draw (4.5,.75) node [anchor=north] {$b^\vee$};
  \draw (5,0) -- (5,.25);
  \draw[<-] (4,0) -- (4,.25);

  \draw[->] [shift={+(7,0)}](0,0) .. controls (1,1) .. (0,2);
  \draw[<-] [shift={+(7,0)}](1,0) .. controls (0,1) .. (1,2);
  \draw [shift={+(7,0)}](1.5,1) node {=};
  \draw[->] [shift={+(7,0)}](2,0) --(2,2);
  \draw[<-] [shift={+(7,0)}](3,0) -- (3,2);
  \end{tikzpicture}
\end{equation}

\begin{equation}\label{eq:rel3}
  \begin{tikzpicture}[>=stealth,baseline=-5pt]
  \draw [shift={+(0,0)}](0,0) arc (180:360:0.5cm);
  \draw [shift={+(0,0)}][->](1,0) arc (0:180:0.5cm);
  \filldraw [shift={+(1,0)}][blue](0,0) circle (2pt);
  \draw [shift={+(0,0)}](1,0) node [anchor=east] {$b$};
  \draw [shift={+(0,0)}](2,0) node{$= \tr(b)$};

  \draw  [shift={+(5,0)}](0,0) .. controls (0,.5) and (.7,.5) .. (.9,0);
  \draw  [shift={+(5,0)}](0,0) .. controls (0,-.5) and (.7,-.5) .. (.9,0);
  \draw  [shift={+(5,0)}](1,-1) .. controls (1,-.5) .. (.9,0);
  \draw[->] [shift={+(5,0)}](.9,0) .. controls (1,.5) .. (1,1);
  \draw  [shift={+(5,0)}](1.5,0) node {$=$};
  \draw  [shift={+(5,0)}](2,0) node {$0$};
  \end{tikzpicture}
\end{equation}

In the first equation in~\eqref{eq:rel2}, the summation is taken over a basis $\mathcal{B}$ of $B^\G$, and this morphism is easily seen to be independent of the choice of basis.  We assign a $\Z$-grading to the space of planar diagrams by defining

\[
  \begin{tikzpicture}[>=stealth]
  \draw  (-.5,.5) node {$\deg$};
  \draw [->](0,0) -- (1,1);
  \draw [->](1,0) -- (0,1);
  \draw (1.5,.5) node{$ = 0,$};
  \end{tikzpicture}
\]
\[
  \begin{tikzpicture}[>=stealth]
  \draw  (-.5,-.25) node {$\deg$};
  \draw[->] (0,0) arc (180:360:.5);
  \draw (1.75,-.25) node{$ = \deg$};
  \draw[->] (3.5,-.5) arc (0:180:.5);
  \draw (4.5,-.25) node{$ =-1,$};
  \end{tikzpicture}
\]
\[
  \begin{tikzpicture}[>=stealth]
  \draw  (-.5,-.25) node {$\deg$};
  \draw[<-] (0,0) arc (180:360:.5);
  \draw (1.75,-.25) node{$ = \deg$};
  \draw[<-] (3.5,-.5) arc (0:180:.5);
  \draw (4.5,-.25) node{$ = 1,$};
  \end{tikzpicture}
\]
and by defining the degree of a dot labeled by $b$ to be the degree of $b$ in the graded algebra $B^\G$.  When equipped with these assignments, all of the graphical relations are homogeneous, and composition of morphisms is compatible with the grading.

Just as in Section~\ref{subsec:graphical-h-cat}, we denote by $\cH_\G$ the Karoubi envelope of $\cH'_\G$ (see Definition~\ref{def:Karoubi-envelope}).  Since $\cH_\G$ is a graded category, the (split) Grothendieck group $K_0(\cH_\G)$ of $\cH_\G$ is an algebra over $\Z[t,t^{-1}]$, where multiplication by $t$ corresponds to the shift $\langle 1 \rangle$.

It follows from the local relations in $\cH'_\G$ that upward oriented crossings satisfy the symmetric group relations, while dots labeled by elements of $\G$ satisfy the relations in $\G$.  Furthermore, upward crossings and dots labeled by elements of $\G$ satisfy the relations of the wreath product $S_n\wr \G :=
S_n \rtimes \G^n$.
Thus we have a canonical homomorphism
\begin{equation} \label{eq:S_n-to-Q+}
  \C[S_n\wr \G] \to \End_{\cH'_\G} (Q_{+^n}).
\end{equation}
Similarly, since each space of morphisms in $\cH'_\G$ consists of diagrams up to isotopy, downward oriented crossings and dots labeled by elements of $\G$ also satisfy the wreath product relations and give us a canonical homomorphism
\begin{equation} \label{eq:S_n-to-Q-}
  \C[S_n\wr \G] \to \End_{\cH'_\G} (Q_{-^n}).
\end{equation}
As explained in \cite[\S3.1.1]{CauLic10}, to each irreducible representation $i \in I_\G$ of $\G$, there is a naturally associated idempotent $e_i(n)\in \C[S_n\wr \G]$.  We will use the notation $e_i(n)$ to
also denote the images of this idempotent in $\End_{\cH'_\G}
(Q_{+^n})$ and $\End_{\cH'_\G} (Q_{-^n})$ under the canonical
homomorphisms~\eqref{eq:S_n-to-Q+} and~\eqref{eq:S_n-to-Q-}.  We then define the following objects in $\cH'_\G$:
\[
  S_{i,+}^n = (Q_{+^n}, e_i(n)),\quad S_{i,-}^n = (Q_{-^n}, e_i(n)).
\]

\begin{theo}[{\cite[Th.~1]{CauLic10}}] \label{thm:main-object-isom-2}
In the category $\cH_\G$, for $i,j \in I_\Gamma$, we have
\begin{align*}
    S_{i,+}^n \otimes S_{j,+}^m &\cong S_{j,+}^m \otimes S_{i,+}^n, \\
  S_{i,-}^n \otimes S_{j,-}^m &\cong S_{j,-}^m \otimes S_{i,-}^n, \\
  S_{i,-}^n \otimes S_{i,+}^m &\cong \bigoplus_{k\geq 0} S_{i,+}^{m-k} \otimes S_{i,-}^{n-k}\otimes \rmH^*(\mathbb{P}^k),\\
  S_{i,-}^n \otimes S_{j,+}^m &\cong \bigoplus_{k=0,1} S_{j,+}^{m-k} \otimes S_{i,-}^{n-k}, \text{ when }
  \langle i,j\rangle = -1,\\
  S_{i,-}^n \otimes S_{j,+}^m &\cong S_{j,+}^m \otimes S_{i,-}^n \text{ when }   \langle i,j\rangle = 0.
\end{align*}
We thus have a well-defined algebra homomorphism $\bF : \fh^{L_{\G,t}}_\Z \to K_0(\cH_\G)$, given by
\[
  \bF(b_i^{(n)}) = [S_{i,-}^n],\quad \bF(a_i^{(n)}) = [S_{i,+}^n].
\]
The homomorphism $\bF$ is an isomorphism.
\end{theo}
In the above, the cohomology of projective space $\rmH^*(\mathbb{P}^k)$ is notation for a direct sum of copies of shifts of the identity, symmetric about the origin, so that, for an object $A$ of $\cH_\G$, we have
\[
  A \otimes 	\rmH^*(\mathbb{P}^k) \cong A\langle-k\rangle \oplus A\langle-k+2\rangle \oplus \dots \oplus A\langle k-2\rangle \oplus A\langle k\rangle.
\]
The shifts $\langle k \rangle$ should not be confused with the pairing $\langle i,j \rangle$ on the lattice.

Note that in Theorem~\ref{thm:main-object-isom}, the analogous homomorphism $\bF$ is only known to be injective; that it is an isomorphism is a conjecture.  The reason that $\bF$ can be shown to be an isomorphism in Theorem~\ref{thm:main-object-isom-2} is that each endomorphism algebra in $\cH_\G$ has a nontrivial grading, inherited from the nontrivial grading on the algebra $B^\G$, with finite-dimensional graded pieces.  This allows one to show that endomorphism algebras in $\cH_\G$ have the Krull-Schmidt property, and hence that objects of $\cH_\G$ decompose uniquely as direct sums of indecomposable objects.  The category $\cH(q)$ lacks such a grading, which makes identifying the Grothendieck group more difficult.

%%%%%%%%%%%%%%%%%%%%%%%%%%%%%%%%%%%%%%%%%%%%%%%%%%%%%
\subsection{Categorification of quantized Fock space}
%%%%%%%%%%%%%%%%%%%%%%%%%%%%%%%%%%%%%%%%%%%%%%%%%%%%%

Let $X_\Gamma$ denote the minimal resolution of the singular variety $\C^2/\Gamma$.
By the Nakajima-Grojnowski constructions described in Section~\ref{subsec:coh}, the Heisenberg algebra $\fh^{\rmH^2(X_\Gamma,\Z)}$ acts on the cohomology of all Hilbert schemes, giving a geometrization of the Fock space representation,
\[
  \cF^{\rmH^2(X_\Gamma,\Z)} \cong \bigoplus_{n=0}^\infty \rmH^*(X_\Gamma^{[n]},\C) \qquad \text{(as $\fh^{\rmH^2(X_\Gamma,\Z)}$-modules)}.
\]

The multiplicative group $\C^*$ acts on $\C^2$ by scaling, and this action commutes with the action of $\Gamma$.  Thus $X_\Gamma^{[n]}$ inherits an action of $\C^*$.
In view of this fact, and by analogy with geometric constructions of representations of quantum affine algebras via quiver varieties \cite{Nak01}, it is natural to suspect that the quantized Heisenberg algebra $\fh^{L_{\G,t}}$ acts on the $\C^*$-equivariant $K$-theory of these Hilbert schemes.  In fact, a much stronger statement is true, as we now explain.

Let $D^b_{\C^*} (X_\G^{[n]})$ denote the bounded derived category of $\C^*$-equivariant coherent sheaves on $X_\G^{[n]}$.  For $n \in \N$, let $\mathcal{D}_n = \bigoplus_{m\geq 0} D^b_{\C^*} (X_\G^{[n]}\times X_\G^{[m]})$.  In \cite[\S4]{CauLic10}, a functor $\bA^\G_n : \cH_\G \to \mathcal{D}_n$ is defined.  (Note that there is some inconsistency between the notation here and the notation in \cite{CauLic10}; in that paper, the notation $A^\G_n$ was used for an algebra, and notation for the functor from $\cH_\Gamma$ to $\mathcal{D}_n$ was not introduced.)

Now, for $m,n \in \N$, there is a natural functor from $D^b_{\C^*} (X_\G^{[n]}\times X_\G^{[m]})$ to $\Fun(D^b(X_\G^{[n]}),D^b(X_\G^{[m]}))$.  If $A \in D^b_{\C^*} (X_\G^{[n]}\times X_\G^{[m]})$, then the associated functor, which is known as the \emph{Fourier-Mukai transform}, is defined by
\[
	B \mapsto {\pi_2}_*(\pi_1^*(B) \otimes A),
\]	
where $\pi_1$ and $\pi_2$ are the projections from $X_\G^{[n]}\times X_\G^{[m]}$ to the first and second factors respectively, and all operations are derived.  In general, when the spaces involved are not compact, one needs to take care when defining the functor associated to a kernel.  However, in this case the objects of interest in  $\mathcal{D}_n$ are proper over both factors, and the non-compactness of $X_\G$ does not cause any technical difficulty.  We refer the reader to \cite{H06} for more information about derived categories of coherent sheaves and Fourier-Mukai transforms.

We now define the functor $\bA_\G$ to be the composition
\begin{multline}
  \bA_\G : \cH_\G \xrightarrow{\bigoplus_{n \in \N} \bA_n^\Gamma} \bigoplus_{m,n \in \N} {D^b_{\C^*} (X_\G^{[n]}\times X_\G^{[m]})} \to \bigoplus_{m,n \in \N} \Fun(D^b(X_\G^{[n]}),D^b(X_\G^{[m]})) \\
  \xrightarrow{} \Fun \left( \bigoplus_{n \in \N} D^b_{\C^*} (X_\G^{[n]}), \bigoplus_{n \in \N} D^b_{\C^*} (X_\G^{[n]}) \right).
\end{multline}

Thus $\bA_\Gamma$ defines a representation of the category $\cH_\G$ on the category $\bigoplus_n D^b_{\C^*} (X_\G^{[n]})$.

\begin{theo}[Categorification of quantized Fock space]
We have a commutative diagram
\begin{equation} \label{eq:Fock-space-cat-diagram-2}
  \xymatrix{
    \cH_\G \ar[r]^(0.18){\bA_\Gamma} \ar[d]_{K_0} &  \Fun \left( \bigoplus_{n \in \N} D^b_{\C^*} (X_\G^{[n]}), \bigoplus_{n \in \N} D^b_{\C^*} (X_\G^{[n]}) \right) \ar[d]^{K_0} \\
    K_0(\cH_\G) \ar[r]^{K_0(\bA_\Gamma)} & \End \mathcal{F}^{L_\G}_\Z \\
    \fh^{L_{\G,t}}_\Z \ar@{->}[u]^{\bF}_\cong \ar[ur]
  },
\end{equation}
where the arrow $\fh^{L_{\G,t}}_\Z \to \End \mathcal{F}^{L_\G}_\Z$ is the action of the integral form of the quantized Heisenberg algebra on the integral form of the quantized Fock space.
\end{theo}

As an immediate corollary, we recover an action of the quantized Heisenberg algebra $\fh^{L_{\G,t}}$ on the $\C^*$-equivariant $K$-theory of Hilbert schemes, answering a question of Nakajima \cite[Question~8.35]{Nak99}.  Thus the action of the quantized Heisenberg algebra on $\C^*$-equivariant K-theory yields a geometric realization of a quantized Fock space.

The categories $D^b_{\C^*}(X_\G^{[n]})$ can be replaced by equivalent derived categories of modules over a finite-dimensional Koszul algebra, see \cite[\S8]{CauLic10}. After this replacement, the representation of $\cH_\G$ becomes abelian rather than triangulated: the endofunctors assigned to each object are given by explicit flat bimodules (rather than complexes thereof) and the natural transformations are bimodule maps.  This abelian construction, while less geometric than the construction on Hilbert schemes, is closer in spirit to the construction of Section~\ref{subsec:CatFock}.

As with the Heisenberg categorifications of \cite{LS11} and \cite{Kho10}, the category $\cH_\G$ has a considerable amount of structure at the level of morphisms.
For example, although left curls on an upward pointing strand are zero, right curls (which have degree $2$) are interesting morphisms.  As shorthand, we will draw right curls in $\cH_\G$ as hollow dots (see~\eqref{eq:right-curl-shorthand}).  These hollow dots satisfy an ``affine Hecke" type relation with crossings which involves the creation of labeled solid dots,
\[
  \begin{tikzpicture}[>=stealth,baseline=10pt]
  \draw [->](0,0) -- (1,1);
  \draw [->](1,0) -- (0,1);
  \draw [red](.25,.25) circle (4pt);
  \draw (1.25,.5) node{$=$};
  \draw [->](1.5,0) -- (2.5,1);
  \draw [->](2.5,0) -- (1.5,1);
  \draw [red](2.25,.75) circle (4pt);

  \draw (3.2,.5) node{$+\sum_{b \in \mathcal{B}}$};

  \draw [shift={+(4.0,0)}][->](0,0) -- (0,1);
  \draw [shift={+(4.0,0)}][->](.75,0) -- (.75,1);
  \filldraw [shift={+(4.0,0)}][blue](0,.33) circle (2pt);
  \draw [shift={+(4.0,0)}](0,.33) node [anchor=west]{$b$};
  \filldraw [shift={+(4.0,0)}][blue](.75,.66) circle (2pt);
  \draw [shift={+(4.0,0)}] (.75,.66) node [anchor=west]{$b^\vee$};
  \end{tikzpicture},
\]
\[
  \begin{tikzpicture}[>=stealth,baseline=10pt]
  \draw [->](0,0) -- (1,1);
  \draw [->](1,0) -- (0,1);
  \draw [red](.25,.75) circle (4pt);
  \draw (1.25,.5) node{$=$};
  \draw [->](1.5,0) -- (2.5,1);
  \draw [->](2.5,0) -- (1.5,1);
  \draw [red](2.25,.25) circle (4pt);

  \draw (3.2,.5) node{$+\sum_{b \in \mathcal{B}}$};

  \draw [shift={+(4.0,0)}][->](0,0) -- (0,1);
  \draw [shift={+(4.0,0)}][->](.75,0) -- (.75,1);
  \filldraw [shift={+(4.0,0)}][blue](0,.33) circle (2pt);
  \draw [shift={+(4.0,0)}](0,.33) node [anchor=west]{$b$};
  \filldraw [shift={+(4.0,0)}][blue](.75,.66) circle (2pt);
  \draw [shift={+(4.0,0)}] (.75,.66) node [anchor=west]{$b^\vee$};
  \end{tikzpicture},
\]
where the summations are over a basis $\mathcal{B}$ of the finite-dimensional Frobenius algebra $B^\G$.  These relations are reminiscent of relations in the degenerate affine Hecke algebra associated to wreath products, see \cite[\S3.5.1]{CauLic10} and references therein.

%%%%%%%%%%%%%%%%%%%%%%%%%%%%%%%%%%%%%%%%%%%%%%%%%%%%%%%%%%%%%%%%%%%
%
\section{Further historical remarks}
\label{sec:further}
%
%%%%%%%%%%%%%%%%%%%%%%%%%%%%%%%%%%%%%%%%%%%%%%%%%%%%%%%%%%%%%%%%%%%

We have not been able to treat all aspects of the categorification of the Heisenberg algebra in this paper.  We therefore conclude by indicating some other related results that have appeared it the literature.  It would be interesting to further elucidate the connection between the results below and the strong categorifications described above.

In \cite{SV11}, Shan and Vasserot defined an action of the Heisenberg algebra on the Grothen\-dieck group of a certain category of modules for cyclotomic rational double affine Hecke algebras. In this way they obtain a categorification of the Fock space representations.  They then describe the relationship between a certain filtration on the Grothendieck group (by support) and a representation theoretic grading on the Fock space, allowing them to prove a conjecture of Etingof and compute the number of finite-dimensional simple objects in the representation category.

We also hope that the constructions of \cite{CauLic10} can be modified to cover the case where
$\G$ is trivial and $X_\G = \C^2$.  This case is particularly interesting because of the relationship between Hilbert schemes of points on $\C^2$ and many other algebraic structures, including elliptic Hall algebras, shuffle algebras, and Macdonald polynomials.  At the level of localized equivariant K-theory, the Hilbert scheme of points on $\C^2$ has been studied by Feigen-Tsymbaliuk \cite{FT09} and Schiffmann-Vasserot \cite{ScV11}.  It would be interesting to lift their constructions to strong categorifications.

Many other fundamental structures related to Heisenberg algebras and vertex operator algebras have been studied at the level of geometrization and weak categorification. We cannot do justice in this expository paper to the large body of work done by many mathematicians over the last 15 years on this subject.  As an imperfect compromise, we refer the readers to Frenkel-Jing-Wang \cite{FJW2,FJW1}, Carlsson-Okounkov \cite{CO08}, Carlsson \cite{C09}, Lehn-Soerger \cite{LSo01}, Qin-Wang \cite{QW02}, Li-Qin-Wang \cite{LQW05}, Okounkov-Pandharipande \cite{OP10}, and Licata-Savage \cite{LS10}, each of which describes some aspect of the relationship between algebraic structures like Heisenberg algebras and vertex operators and the geometry of Hilbert schemes or representation theory of symmetric groups.  Lifting the mathematics studied in these and many other closely related papers to strong categorifications is an important area of current activity.

%%%%%%%%%%%%%%%%%%%%%%%%%%%%%%%%%%%%%%%%%%%%%%%%%%%%%%%%%%%%%%%%%%%%%
%%%%%%%%%%%%%%%%%%%%%%%%%%%%%%%%%%%%%%%%%%%%%%%%%%%%%%%%%%%%%%%%%%%%%

\bibliographystyle{plain}
\bibliography{heis-biblist}

\end{document}